\documentclass[12pt,reqno]{amsart}
\usepackage{amsmath}
\usepackage[latin1]{inputenc}
\usepackage{amssymb}
\usepackage{amscd}
\usepackage{natbib}
\usepackage{mathrsfs}
\usepackage[pdftex]{graphicx}
\usepackage{graphicx}
\usepackage[OT2, T1]{fontenc}
\usepackage[russian,english]{babel}
\usepackage[babel=true]{csquotes}
\usepackage{chngcntr}
\counterwithout{figure}{section}
\counterwithout{figure}{subsection}

\newcommand{\E}[0]{\ensuremath{\mathbf{E}}}
\newcommand{\F}[0]{\ensuremath{\mathcal{F}}}

\newcommand{\Pe}[0]{\ensuremath{\mathbf{P}}}
\newcommand{\R}[0]{\ensuremath{\mathbb{R}}}

\newcommand{\Z}[0]{\ensuremath{\mathbb{Z}}}

\newtheorem{assu}{\sc{Assumption}}[section]
\newtheorem{theo}{\sc{Theorem}}[section]
\newtheorem{lemm}[theo]{\sc{Lemma}}
\newtheorem{defi}[theo]{\sc{Definition}}

\newtheorem{cor}[theo]{\sc{Corollary}}
\newtheorem{rmq}[theo]{\sc{Remark}}

\author{Alexandre Boritchev}
\title[]{Exponential convergence to the stationary measure for a class of 1D Lagrangian systems with random forcing}
\date{\today}
\begin{document}


\maketitle

\textbf{Abstract.} We prove exponential convergence to the stationary measure for a class of 1d Lagrangian systems with random forcing in the space-periodic setting:
\begin{equation} \nonumber
\phi_t+\phi_x^2/2=F^{\omega},\ x \in S^1=\R/\Z.
\end{equation}
This confirms a part of a conjecture formulated in \cite{GIKP05}. Our result is a consequence (and the natural stochastic PDE counterpart) of the results obtained in \cite{BK13, EKMS00}. It is also the natural analogue of the deterministic result \cite{ISM09} which holds in a generic setting.

\tableofcontents

\section*{Abbreviations}

\begin{itemize}
\item 1d: one-dimensional
\item a.s.: almost surely
\item i.i.d.: independent identically distributed
\item r.v.: random variable
\end{itemize}

\section*{Introduction} \label{intro}

We are concerned with 1d random Lagrangian systems of the mechanical type, i.e. of the form:
\begin{equation} \label{L}
L^{\omega}(x,v,t)=v^2/2+F^{\omega}(x,t),\ x \in S^1=\R/\Z,
\end{equation}
where $F^{\omega}(x,t)$ is a smooth function in $x$ and a stationary random process in $t$ (of the kick or white force type: see Section~\ref{random}).
The Legendre-Fenchel transform gives us the corresponding Hamiltonian:
\begin{equation} \label{H}
H^{\omega}(x,p,t)=p^2/2-F^{\omega}(x,t),
\end{equation}
and the Hamilton-Jacobi equation:
\begin{eqnarray} \label{HJ}
&\phi_t + \phi_x^2/2= F^{\omega}.
\end{eqnarray}
Here, we consider only the case where $\phi$ is itself $1$-periodic in space. In this case the function $u=\phi_x$ satisfies the randomly forced inviscid Burgers equation:
\begin{eqnarray} \label{Bur}
&u_t + uu_x=  (F^{\omega})_x,\ x \in S^1=\R/\Z.
\end{eqnarray}
Note that it is equivalent to consider a solution of (\ref{Bur}) and a solution of (\ref{HJ}) defined up to an additive constant. Under the assumptions of Section~\ref{random}, both of these equations are well-posed  and their solutions define Markov processes. Therefore, we can consider the corresponding stationary measure. Its existence and uniqueness has been proved by E, Khanin, Mazel and Sinai in the white force case in the seminal work \cite{EKMS00} using the Lagrangian representation of the solutions. For the case of a kick force, see for instance \cite{IK03}. This result was clarified and generalised to the multi-d case by Khanin and his collaborators in \cite{BK13, GIKP05, KZ11} with more transparent assumptions on the forcing. In these papers, no explicit estimates on the speed of convergence are given. However, an exponential bound on the speed of convergence locally in space far from the shocks has been obtained in \cite{BFK00}. In the papers mentioned above, the key object is the global minimiser and the key fact is its hyperbolicity.
\\ \indent
We have previously obtained a bound on the speed of convergence to the stationary measure for solutions of (\ref{Bur}) both in the 1d and in the multi-d case in \cite{BorK,BorW, BorM} using stochastic PDE techniques. This bound is of the type $C t^{-\delta/p},\ \delta>0$, in the dual-Lipschitz metric corresponding to the space of functions in $L_{p},\ p \in [1,\infty)$. Although this bound is given for the equation with an additional viscous term $\nu u_{xx}$, it is independent from the viscosity coefficient $\nu$, and thus it still holds when we pass to the limit $\nu \rightarrow 0$.
\\ \indent
In our paper, we give an exponential bound for the speed of convergence to the stationary measure for solutions of (\ref{Bur}) in the dual-Lipschitz metric mentioned above, which gives a partial answer to the conjecture stated in \cite[Section 4]{GIKP05}. The part of the conjecture which remains open is proving that this exponential bound still holds if we add a positive viscosity coefficient $\nu$, uniformly in $\nu$. The main technical difficulty is that introducing this term destroys the well-understood structure of the minimisers.
\\ \indent
It is very likely that the estimate we obtain is optimal since it gives exactly the same results as the optimal bound obtained in the generic nonrandom case by R.~Iturriaga and H.~Sanchez-Morgado \cite{ISM09}. Note that the metrics we use are also optimal since it is impossible to obtain such a bound in the Lipschitz-dual space corresponding to $L_{\infty}$ for solutions of (\ref{Bur}) which are discontinuous with a strictly positive probability.
\\ \indent
Finally, we would like to emphasize the link between our work and the corresponding deterministic results belonging to the realm of the weak KAM theory developed by A.~Fathi and J.~Mather \cite{Fat05}. In particular, there is a striking correspondence between the scheme of our proof and the one in \cite{ISM09}, which follows a general rule: the results which hold in the random case under fairly weak assumptions are similar to the results which hold in the nonrandom case under more stringent genericity assumptions. For more on this subject and the link with the Aubry-Mather theory, see for instance \cite{IK03}.

\begin{rmq}
Our results extend to the case where $\phi$, instead of being periodic in space, satisfies
$$
\phi(x+1)=\phi(x)+b,\ x \in \R.
$$
The proofs are exactly the same since we use the results of \cite{BK13, EKMS00}, which hold for all values of $b$.
\\ \indent
Our results extend to a class of non-mechanical convex in $p$ Hamiltonians of the type $H(p)+F^{\omega}(t,x)$ with $F^{\omega}$ as above, under some assumptions of the Tonelli type \cite{Fat05}.
\\ \indent
\end{rmq}

\section{Notation and setting} \label{nota}

\subsection{Random setting} \label{random}

We consider the mechanical Hamilton-Jacobi equation with two different types of additive forcing in the right-hand side and a $C^{\infty}$-smooth initial condition $\phi^0$.
\\ \indent
We begin by formulating the assumptions on potentials, which are (except \ref{C1} (i)) the same as in the paper \cite{BK13}:

\begin{assu} \label{C1}
In the \enquote{kicked} case, we assume the following.
\\
(i)\ The kicks at integer times $j$ are of the form
$$
F^{\omega}(j)(x)=\sum_{k=1}^{K}{c_k^{\omega}(j) F^k(x)},
$$
where $F^k$ are $C^{\infty}$-smooth potentials on $S^1=\R/\Z$. The random vectors
\\
$(c_k^{\omega}(j))_{1 \leq k \leq K}$ are i.i.d. $\R^K$-valued random variables defined on a probability space $(\Omega, \F, \Pe)$. Their distribution on $\R^K$, denoted by $\lambda$, is assumed to be absolutely continuous with respect to the Lebesgue measure, and all of its moments are assumed to be finite.
\\
(ii)\ 
The potential $0$ belongs to the support of $\lambda$.
\\
(iii)\ The mapping from $S^1$ to $\R^K$ defined by
$$
x \mapsto (F^1(x),...,F^K(x))
$$
is an embedding.
\end{assu}

\begin{assu} \label{C2}
In the case of the white force potential, we assume the following.
\\
(i)\ The forcing has the form
$$
F^{\omega}(x,t)=\sum_{k=1}^{K}{(W_k^{\omega})_t(t) F^k(x)},
$$
where $F^k$ are smooth potentials on $S^1$, and $(W^{\omega}_k)_t$ are independent white noises defined on a probability space $(\Omega, \F, \Pe)$, i.e. weak time derivatives of independent Wiener processes $W_k^{\omega}(t)$.
\\
(ii)\ The mapping from $S^1$ to $\R^K$ defined by
$$
x \mapsto (F^1(x),...,F^K(x))
$$
is an embedding.
\end{assu}

\begin{rmq}
Our results also extend to the case of infinite-dimensional noise: here the necessary restriction is that the noise remains smooth in space. For instance, we can put independent white noises on each Fourier mode in such a way that the amplitude of the noise decreases exponentially with the wavenumber.
\end{rmq}
\smallskip
\indent
In the white noise case, we denote by $G$ an antiderivative in time of the forcing:
$$
G^{\omega}(x,t)=\sum_{k=1}^{K}{W_k^{\omega}(t) F^k(x)},
$$
where $W_k^{\omega}(t)$ are independent standard Wiener processes with $W_k^{\omega}(0)=0$. Since we will only consider time differences of $G$, the particular choice of antiderivative has no importance.
\\ \indent
In both cases, $F^{\omega}$ will be abbreviated as $F$, and in the white force case $F(\cdot,t)$ will be abbreviated as $F(t)$, and similarly for $G$.
\\ \indent
Note that since we have $G(t) \in C^{\infty}$ for every $t$, a.s., we can redefine the forcing $F$ in the white force case so that this property holds for all $\omega \in \Omega$.

\subsection{Functional spaces and Sobolev norms} \label{Sob}

Consider an
\\
integrable function $v$ on $S^1$. For $p \in [1,\infty]$, we denote its $L_p$ norm by $\left|v\right|_p$. The $L_2$ norm is denoted by  $\left|v\right|$, and $\left\langle \cdot,\cdot\right\rangle$ stands for the $L_2$ scalar product. From now on $L_p,\ p \in [1,\infty],$ denotes the space of zero mean value functions in $L_p(S^1)$. Similarly, $C^{\infty}$ is the space of $C^{\infty}$-smooth zero mean value functions on $S^1$.
\\ \indent
For a nonnegative integer $m$ and $p \in [1,\infty]$, $W^{m,p}$ stands for the Sobolev space of zero mean value functions $v$ on $S^1$ with finite homogeneous norm
\begin{equation} \nonumber
\left|v\right|_{m,p}=\left|\frac{d^m v}{dx^m}\right|_p.
\end{equation}
In particular, $W^{0,p}=L_p$ for $p \in [1,\infty]$. We will never use Sobolev norms for non-zero mean functions: in particular, for solutions of (\ref{HJ}) we will only consider the Lebesgue norms.
\\ \indent
Since the length of $S^1$ is $1$, we have:
$$
|v|_1 \leq |v|_{\infty} \leq |v|_{1,1} \leq |v|_{1,\infty} \leq \dots \leq |v|_{m,1} \leq |v|_{m,\infty} \leq \dots
$$
\indent
We denote by $L_{\infty}/\R$ the space of functions in $L_{\infty}$ defined modulo an additive constant endowed with the norm:
$$
|u-v|_{L_{\infty}/\R}=\inf_{K \in \R}{|u-v-K|_{\infty}}
$$
\\ \indent
We recall a version of the classical Gagliardo--Nirenberg inequality (see \cite[Appendix]{DG95}):

\begin{lemm} \label{GN}
For a smooth zero mean value function $v$ on $S^1$,
$$
\left|v\right|_{\beta,r} \leq C \left|v\right|^{\theta}_{m,p} \left|v\right|^{1-\theta}_{q},
$$
where $m>\beta\geq 0$, and $r$ is defined by
$$
\frac{1}{r}=\beta-\theta \Big( m-\frac{1}{p} \Big)+(1-\theta)\frac{1}{q},
$$
under the assumption $\theta=\beta/m$ if $p=1$ or $p=\infty$, and $\beta/m \leq \theta < 1$ otherwise. The constant $C$ depends on $m,p,q,\beta,\theta$.
\end{lemm}

Subindices $t$ and $x$, which can be repeated, denote partial differentiation with respect to the corresponding variables. We denote by $v^{(m)}$ the $m$-th derivative of $v$ in the variable $x$. For brevity, the function $v(t,\cdot)$ is denoted by $v(t)$.

\subsection{Agreements} \label{agree}

All functions which we consider in this paper are real-valued. All quantities denoted by $K_i$, $i$ being a natural number, are positive constants which only depend on the general features of the system (i.e. the statistical distribution of the forcing): they are nonrandom and do not depend on the initial condition. All quantities denoted by $C_i$, $i$ being a natural number, are (time-dependent) r.v.'s with all moments finite and uniformly bounded in time, which do not depend on the initial condition.
\\ \indent
We will always denote by $\phi(t,x)$ a solution of (\ref{HJ}) and by $u(t,x)$ its derivative, which solves (\ref{Bur}), respectively for initial conditions $\phi^0$ and $u^0=\phi^0_x$. We will denote accordingly the solutions for two initial conditions $\phi^0,\overline{\phi^0}$. The assumptions on the forcing are the ones given in Section~\ref{random}.

\section{Dynamical objects and stationary measure} \label{dyn}

Here we introduce the Lagrangian dynamical objects in the setting described in the previous section. Note that all the results of Sections~\ref{min}-\ref{meas} hold under much more general assumptions: for instance, it is possible to drop (iii) in Assumption~\ref{C1} or (ii) in Assumption~\ref{C2}. However, these hypotheses will be extremely important for the results which will be given in Section~\ref{scheme}. For more details on the definitions given below, see \cite{GIKP05, IK03}.

\subsection{Lagrangian formulation and minimisers} \label{Lagrange}

\begin{defi} \label{min}
For a time interval $[s,t]$, we say that a curve $\gamma_{s,t}^{y,x}(\tau)$ is a \textbf{minimiser} if it minimises the action
\begin{align} \nonumber
&A(\gamma)=\frac{1}{2} \int\limits_{s}^{t}{\gamma_t(\tau)^2 d \tau} + \sum_{n \in (s,t])} {\Big(F^n(\gamma(n))\Big)}
\end{align}
in the \enquote{kicked} case and the action
\begin{align} \nonumber
A(\gamma) =& \frac{1}{2} \int\limits_{s}^{t}{\gamma_t(\tau)^2 d \tau}+\int\limits_{s}^{t} { \Bigg( \gamma_t(\tau) \Big(\frac{\partial G}{\partial x}(\gamma(\tau),s)-\frac{\partial G}{\partial x}(\gamma(\tau),\tau) \Big) \Bigg) d \tau }
\\ \nonumber
&+ \Big( G(\gamma(t),t)-G(\gamma(t),s) \Big)
\end{align}
in the white force case, respectively, over all absolutely continuous
\\
curves $\gamma$ such that $\gamma(t)=x$ and $\gamma(s)=y$.
\end{defi}

\begin{defi} \label{phimin}
For a time interval $[s,t]$ and a continuous function $\phi: S^1 \rightarrow \R$, we say that a curve $\gamma_{s,t,\phi}^{x}(\tau): [s,t] \rightarrow S^1$ is a \textbf{$\phi$-minimiser} if it minimises 
$A(\gamma)+\phi(\gamma(s))$
over all absolutely continuous curves on $[s,t]$ such that $\gamma(t)=x$. In particular, all $\phi$-minimisers are minimisers.
\end{defi}

\indent
Now we can define the (pathwise) solution to (\ref{HJ}) for a given $\omega \in \Omega$ and a given continuous initial condition.
Note that by a compactness argument, one can show that for any given endpoint $x$, a minimiser $\gamma$ on $[s,t]$ such that $\gamma(t)=x$ exists. In the white force case, this minimiser gives a time-continuous solution in $L_1$, whereas in the \enquote{kicked} case the solution is a cadlag in time (right-continuous and with a limit to the left) $L_1$-valued function.

\begin{defi} \label{sol}
For a time interval $[s,t]$ and a continuous initial condition $\phi(s): S^1 \rightarrow \R$, for every $\omega$ by definition the (pathwise) solution $\phi: [s,t] \times S^1 \rightarrow \R$ of (\ref{HJ}) is defined by the $\omega$-depended action $A$:
$$
\phi(\tau,x)=A(\gamma)+\phi(s,\gamma(s)),\ \tau \in [s,t],
$$
where $\gamma$ is an $\omega$-dependent $\phi(s)$-minimiser defined on $[s,\tau]$ satisfying $\gamma(\tau)=x$.
\end{defi}

It is easy to check that this solution will verify the semigroup property: in other words, one can define a solution operator
$$
S_{t_1}^{t_2}:\ \phi(t_1) \mapsto \phi(t_2),\ s \leq t_1 \leq t_2 \leq t,
$$
such that for $t_1 \leq t_2 \leq t_3$,
$$
S_{t_2}^{t_3} \circ S_{t_1}^{t_2} = S_{t_1}^{t_3}.
$$
In particular, the following holds:

\begin{lemm} \label{semigroup}
For any $\tau \in (s,t)$, the restriction of any $\phi(s)$-minimiser defined on $[s,t]$ on the time interval $[\tau,t]$ is a $S_{s}^{\tau} \phi(s)$-minimiser.
\end{lemm}

\begin{rmq}
Note that the solution $\phi$ is the limit in $L_1$ of the strong solutions to the equation obtained if we add a viscous term $\nu \phi_{xx}$ to (\ref{HJ}) and then we make $\nu$ tend to $0$ (see \cite{GIKP05}).
\end{rmq}

\begin{defi} \label{onesidedmin}
For a time $t$ and a point $x \in S^1$, we say that a curve $\gamma_{t}^{x,+}(\tau): [t,+\infty) \mapsto S^1$ is a forward one-sided minimiser if it minimises 
$A(\gamma)$
over all absolutely continuous curves such that $\gamma(t)=x$ for compact in time perturbations.
\\ \indent
Namely, we require that if for a curve $\tilde{\gamma}$ such that $\tilde{\gamma}(t)=x$ there exists $T$ such that $\tilde{\gamma}(s) \equiv \gamma(s)$ for $s \geq T$, then $A(\gamma)-A(\tilde{\gamma}) \leq 0$ (this difference is well-defined since it is equal to the difference of the actions on the finite interval $[t,T]$).
%\\ \indent
%Backward one-sided minimisers $\gamma_{t}^{x,-}(\tau): (-\infty,t] \mapsto S^1$ are defined in the same way.
\end{defi}

%\begin{defi} \label{twosidedmin}
%A curve $\gamma(\tau): (-\infty,+\infty) \mapsto S^1$ is a two-sided minimiser if it minimises $A(\gamma)$ up to compact perturbations over all absolutely continuous curves, in the same sense as in Definition~\ref{onesidedmin}.
%\end{defi}

\subsection{Stationary measure and related issues} \label{meas}

Here we give a few results which hold under weak assumptions and are sufficient to ensure that the stationary measure corresponding to (\ref{Bur}) exists and is unique. These results are not new and hold both in the one-dimensional and in the multi-dimensional setting: see \cite{EKMS00,IK03}. Estimates for the speed of convergence are given in \cite{BorK,BorW, BorM}, where all proofs are stated for $\nu >0$, but still hold for $\nu=0$ \cite{EKMS00, GIKP05}. Up to some natural modifications due to the fact that the forcing is now discrete in time, the convergence estimates can be generalised to the kick force case in 1d \cite{BorK}. For more details, see also \cite{KuSh12}, where a random forcing is introduced in a similar setup.
\\ \indent
The flow corresponding to (\ref{Bur}) induces a Markov process, and then we can define the corresponding semigroup denoted by $S_t^{*}$, acting on Borel measures on any $L_p,\ 1 \leq p < \infty$. A \textit{stationary measure} for (\ref{Bur}) is a Borel probability measure defined on $L_p$, invariant with respect to $S_t^{*}$ for every $t$.  A \textit{stationary solution} of (\ref{Bur}) is a random process $v$ defined for $(t,\omega) \in [0,+\infty) \times \Omega$, satisfying (\ref{Bur}) and taking values in $L_p$, such that the distribution of $v(t)$ does not depend on $t$. This distribution is automatically a stationary measure.
\\ \indent
Existence of a stationary measure for (\ref{Bur}) is obtained using uniform bounds for solutions in $W^{1,1}$, which is compactly injected into $L_p,\ p \in [1,\infty)$, and the Bogolyubov-Krylov argument. It is more delicate to obtain uniqueness of a stationary measure, which implies uniqueness for the probability distribution function of a stationary solution.

\begin{rmq}
The most natural space for our model would be the space $L_{\infty}/\R$, on which we could have treated directly the solutions to the equation (\ref{HJ}). Moreover, this is the space in which exponential convergence to the unique stationary solution is proved in the deterministic generic setting in \cite{ISM09}. However, this space is not separable, which makes it delicate to deal with the stationary measure.
\end{rmq}

\begin{defi}
Fix $p \in [1,\infty)$. For a continuous function
$$
g:\ L_p \rightarrow \R,
$$
we define its Lipschitz norm as
$$
|g|_{L(p)}:=|g|_{Lip}+\sup_{L_p}{|g|},
$$
where $|g|_{Lip}$ is the Lipschitz constant of $g$. The set of continuous functions with finite Lipschitz norm will be denoted by $L(p)$.
\end{defi}

\begin{defi}
For two Borel probability measures $\mu_1,\mu_2$ on $L_p$, we denote by $\Vert \mu_1-\mu_2 \Vert^*_{L(p)}$ the Lipschitz-dual distance:
$$
\Vert \mu_1-\mu_2 \Vert^*_{L(p)}:=\sup_{g \in L(p),\ |g|_{L(p)} \leq 1}{\Big| \int_{S^1}{g d \mu_1}-\int_{S^1}{g d \mu_2} \Big|}.
$$
\end{defi}

The following result proved in \cite{BorK,BorW, BorM} is, as far as we are aware, the first explicit estimate for the speed of convergence to the stationary measure of (\ref{Bur}) which is uniform with respect to the viscosity coefficient $\nu$ and is formulated in terms of Lebesgue spaces only. It holds both in the white force and in the kick-force setting in 1d, and only in the white force case in the multidimensional setting. However, a result which holds in the $L_{\infty}$ norm locally in space far from the shocks in 1d has been obtained by Bec, Frisch and Khanin in \cite{BFK00}. The proof in \cite{BorK,BorW, BorM} uses a version of the coupling argument due to Kuksin and Shirikyan \cite[Chapter 3]{KuSh12}. The situation is actually simpler than for the stochastic 2D Navier Stokes equation. In particular, in our setting the "damping time" needed to make the distance between two solutions corresponding to the same forcing small does not depend on the initial conditions. Moreover sincethe flow of (\ref{Bur}) is $L_{1}$-contracting, the coupling argument is simplified.

\begin{theo} \label{algCV}
There exists $\delta>0$ such that for every $p \in [1,\infty)$, there exists a positive constant $K'(p)$ such that we have:
\begin{equation} \\
\Vert S_t^{*} \mu_1- S_t^{*} \mu_2 \Vert^*_{L(p)} \leq K' t^{- \delta/p},\qquad t \geq 1,
\end{equation}
for any probability measures $\mu_1$, $\mu_2$ on $L_p$.
\end{theo}

\subsection{Main results and scheme of the proof} \label{scheme}

Now we are ready to state the main result of the paper.

\begin{theo} \label{main}
There exists $K>0$ such that for every $p \in [1,\infty)$, there is a positive constant $K'(p)$ such that we have:
\begin{equation} \\
\Vert S_t^{*} \mu_1- S_t^{*} \mu_2 \Vert^*_{L(p)} \leq K'(p) exp (- K t/p),\qquad t \geq 0,
\end{equation}
for any probability measures $\mu_1$, $\mu_2$ on $L_p$.
\end{theo}

\indent
The scheme of the proof is, in the spirit, similar to the proof of \cite[Theorem 1]{ISM09}. In that paper the authors use the objects of
the weak KAM theory such as the Peierls barrier, which do not have any directly available counterparts in our setting. However, there is a straightforward dynamical interpretation of their proof in the simplest case. Namely, consider a mechanical Lagrangian
$$
v^2/2+V(x)
$$
such that the potential $V$ is smooth and generic (i.e. it reaches its minimum at a unique point $y$ with $V''(y)>0$).
\\ \indent
The basic idea is that a curve which minimises energy during an interval of time $[0,T]$ remains in a small neighbourhood of $y$ on $[\tau,T-\tau]$ (with $\tau$ $T$-independent). Consequently, since $y$ is a nondegenerate minimum for $V$, we obtain by linearising the Euler-Lagrange equation that at the time $T/2$, all minimisers (independently of the initial condition) are $C \exp(-CT)$-close to $y$, and then a standard argument allows us to conclude that for any initial conditions $\phi^0,\ \overline{\phi^0}$, the solutions of (\ref{HJ}) at time $T$ are $C \exp(-CT)$-close up to an additive constant, i.e.:
$$
\sup_{\phi^0,\ \overline{\phi^0} \in C^0} \inf_{K \in \R} \Big| \phi(T,x)-\overline{\phi}(T,x)-K \Big|_{\infty} \leq C \exp(-CT),
$$
where $K=\phi(T,y)-\overline{\phi}(T,y)$ is $x$-independent.
\\ \indent
In our setting, there are two main ingredients in the proof. Roughly speaking, the first one tells us that for a given initial condition $\phi^0$, the $\phi^0$-minimisers concentrate exponentially. The second one tells us that the one-sided minimisers, which are limits of the $\phi^0$-minimisers on $[0,T]$ as $T \rightarrow +\infty$ for \textit{any} initial condition $\phi^0$, also concentrate exponentially.
\\ \indent
Now we introduce some definitions.
\medskip
\\ \indent
The diameter of a closed set $Z$ can be thought of as the minimal length of a closed interval on $S^1$ containing $Z$.

\begin{defi}
Consider a closed subset $Z$ of $S^1$. Let $a(Z)$ denote the maximal length of a connected component of $S^1-Z$. We define the diameter of $Z$ as
$$
d(Z)=1-a(Z).
$$ 
\end{defi}

\begin{defi} \label{Omega}
For $-\infty<r<s \leq t<+\infty$ and for a fixed function $\phi^0: S^1 \rightarrow \R$, let $\Omega_{r,s,t,\phi^0}$ be the set of points reached, at the time $s$, by $\phi^0$-minimisers on $[r,t]$:
\begin{eqnarray} \nonumber
&\Omega_{r,s,t,\phi^0}=\lbrace \gamma_{r,t,\phi^0}^x(s),\ x \in S^1\rbrace.
\end{eqnarray}
\end{defi}

Now we give the formulations of the two key lemmas, which are (up to notation and a few simplifications) respectively \cite[Theorem 2.1.]{BK13} and \cite[Lemma 5.6.(a)]{EKMS00}. Although the second lemma is only proved in the white force setting in \cite{EKMS00}, its proof in the kick force setting follows the same lines and is technically simpler.

\begin{lemm} \label{expomega}
There exist a random constant $C_1$ and a constant $K_1$ such that we have the inequality:
\begin{eqnarray} \nonumber
& \sup_{\phi^0 \in C^0} d(\Omega_{0,s,s+s',\phi^0}) \leq C_1 \exp(-K_1 s').
\end{eqnarray}
\end{lemm}

\begin{lemm} \nonumber
There exists a random constant $C_2$ and a constant $K_2$ such that we have:
\begin{equation} \label{expmineq}
\sup_{\tilde{\gamma}_1,\tilde{\gamma}_2 \in \Gamma} |\tilde{\gamma}_1(t)-\tilde{\gamma}_2(t)| \leq C_2 \exp(-K_2 t),\ t \geq 0.
\end{equation}
where $\Gamma$ is the set of all forward one-sided minimisers defined on the time interval $[0,+\infty)$.
\end{lemm}

\begin{cor} \label{expmincor}
Consider an initial condition $\phi^0$ and a time $t>0$. Then there exists a random constant $C_3$ and a constant $K_3$ such that for any $\phi^0$-minimiser $\gamma:[0,2t] \rightarrow S^1$ and any forward one-sided minimiser $\delta:\ [0,+\infty) \rightarrow S^1$  we have:
\begin{equation} \label{expmincoreq}
|\gamma(t)-\delta(t)| \leq C_3 \exp(-K_3 t).
\end{equation}
\end{cor}

\textbf{Proof of Corollary~\ref{expmincor}:} 
As we know from \cite[Section 5]{EKMS00}, extracting a subsequence of minimisers (and in particular of $\phi^0$-minimisers) on $[0,s]$ while letting $s$ go to $+\infty$, one obtains a forward one-sided minimiser. In particular, for every $\epsilon$ there exists $s(\epsilon) \geq 2t$, a $\phi^0$-minimiser $\tilde{\gamma}$ defined on $[0,s]$ and a forward one-sided minimiser $\tilde{\delta}$ on $[0,+\infty)$ such that:
\begin{align} \label{expmincorineqbis}
& |\tilde{\gamma}(t)-\tilde{\delta}(t)| \leq \epsilon.
\end{align}

By  Lemma~\ref{expmineq} we have:
\begin{align} \label{expomegacsq}
& |\delta(t)-\tilde{\delta}(t)| \leq C_1 \exp(-K_1 t),
\end{align}
and by Lemma~\ref{expomega}, since the restriction $\tilde{\gamma}|_{[0,2t]}$ is still a one-sided minimiser, we have:
\begin{align} \label{expmineqcsq}
& |\gamma(t)-\tilde{\gamma}(t)| \leq C_2 \exp(-K_2 t),
\end{align}
\indent
Combining the inequalities (\ref{expmincorineqbis})-(\ref{expmineqcsq}) and using the triangular inequality, and then letting $\epsilon$ go to 0, we get (\ref{expmincoreq}) with $K_3=\min(K_1,\ K_2)$ and $C_3=C_1+C_2$. $\qed$

\section{Proof of Theorem~\ref{main}} \label{proof}

We recall the statement of the main theorem:

\textbf{Theorem~\ref{main}:} There exists $K>0$ such that for every $p \in [1,\infty)$, there exists a positive constant $K'(p)$ such that we have:
\begin{equation} \\
\Vert S_t^{*} \mu_1- S_t^{*} \mu_2 \Vert^*_{L(p)} \leq K'(p) exp (- K t/p),\qquad t \geq 0,
\end{equation}
for any probability measures $\mu_1$, $\mu_2$ on $L_p$.
\bigskip
\\ \indent
The proof of the three auxiliary results given below can be found in \cite{EKMS00}. Alternatively, to prove Lemma~\ref{W11} one can take the $\nu$-uniform estimates in \cite{BorK,BorW,BorM} and consider the limit $\nu \rightarrow 0$.

%\begin{lemm} \label{contract}
%Consider two solutions $\phi$, $\overline{\phi}$ of (\ref{HJ}), corresponding to the same random force but different initial conditions in $L_{\infty}$. For all $t \geq s \geq 0$, we have:
%\begin{equation} \nonumber
%|\phi(t)-\overline{\phi}(t)|_{\infty} \leq |\phi(s)-\overline{\phi}(s)|_{\infty}.
%\end{equation}
%Consequently,
%\begin{equation} \label{contracteq}
%|\phi(t)-\overline{\phi}(t)|_{L_{\infty}/\R} \leq |\phi(s)-\overline{\phi}(s)|_{L_{\infty}/\R}.
%\end{equation}
%\end{lemm}

\begin{lemm} \label{W11}
There is a random constant $C_4$ such that for $t \geq 1$, we have:
$$
\sup_{\phi^0 \in C^0} |u(t)|_{1,1} \leq C_4.
$$
\end{lemm}

\begin{cor} \label{Linfty}
For $t \geq 1$, we have:
$$
\sup_{\phi^0 \in C^0} |u(t)|_{\infty} \leq C_4,
$$
where $C_4$ is the same as above.
\end{cor}

\begin{lemm} \label{gammaprime}
For $t \geq 1$, we have:
$$
\sup_{s \in [t,t+1], \gamma \in \Gamma} |\gamma_t(s)|_{1,1} \leq C_4,
$$
where $\Gamma$ is the set of minimisers defined on $[0,t+1]$, and $C_4$ is the same as above.
\end{lemm}

Moreover, we will need the following lemma, analogous to \cite[Section 3, Fact 1]{EKMS00}.

\begin{lemm} \label{connect}
Consider two minimisers $\gamma_1,\gamma_2$, both defined on $[t,T],\ T \geq t+1,$ and satisfying $\gamma_1(T)=\gamma_2(T)$. There is a random constant $C_5$ such that if for $\epsilon>0$ we have:
$$
|\gamma_1(t)-\gamma_2(t)| \leq \epsilon,
$$
then  we have the following inequality for the actions of the minimisers:
\begin{equation} \nonumber
|A(\gamma_1)-A(\gamma_2)| \leq C_5 (\epsilon+\epsilon^2).
\end{equation}
\end{lemm}

\textbf{Proof:} By symmetry, it suffices to prove that:
\begin{equation} \label{connecteqabs}
A(\gamma_2) \leq A(\gamma_1)+C_5 (\epsilon+\epsilon^2);
\end{equation}
$C_5$ will be fixed later. Consider the curve $\tilde{\gamma}_1: [t,T] \rightarrow S^1$ defined by:
\begin{align} \nonumber
& \tilde{\gamma}_1(s)=\gamma_1(s)+(t+1-s)(\gamma_2(t)-\gamma_1(t)),\ s \in [t,t+1].
\\ \nonumber
& \tilde{\gamma}_1(s)=\gamma_1(s),\ s \in [t+1,T].
\end{align}
Using Definition~\ref{min} and Lemma~\ref{gammaprime}, we get:
$$
A(\tilde{\gamma}_1) \leq A(\gamma_1)+ C_5 (\epsilon+\epsilon^2).
$$
On the other hand, since $\tilde{\gamma}_1$ has the same endpoints as the minimiser $\gamma_2$, we get:
$$
A(\gamma_2) \leq A(\tilde{\gamma}_1).
$$
Combining these two inequalities proves (\ref{connecteqabs}).\ \qed
\\
\bigskip

The proof of the following lemma follows the lines of \cite{ISM09}.

\begin{lemm} \label{expcontract}
Consider two solutions $\phi$ and $\overline{\phi}$ of (\ref{HJ}) defined on the time interval $[0,+\infty)$. There there exist $M>0$ and a random constant $C_6$ such that we have:
$$
\E |\phi(t)-\overline{\phi}(t)|_{L_{\infty}/\R} \leq C_6 \exp(-Mt),\ t \geq 0.
$$
\end{lemm}

\textbf{Proof of Lemma~\ref{expcontract}:} Consider two solutions $\phi$ and $\overline{\phi}$ to (\ref{HJ}) corresponding to the same forcing and different initial conditions at time $0$. Using Definition~\ref{sol}, we get for any $t \geq 1$ and $x \in S^1$:
\begin{align} \label{first}
& \phi(2t,x)-\overline{\phi}(2t,x)
\\ \nonumber
&= \phi(t,\gamma_1(t))+A(\gamma_1|_{[t,2t]})-\overline{\phi}(t,\gamma_2(t))-A(\gamma_2|_{[t,2t]}),
\end{align}
where $\gamma_1$ and $\gamma_2$ are respectively a $\phi^0$- and a $\overline{\phi^0}$-minimiser on $[0,2t]$ ending at $x$. By Corollary~\ref{expmincor}, we have:
\begin{align} \label{expclose}
|\gamma_i(t)-y| \leq C_3 \exp(-K_3 t),\ i=1,2,
\end{align}
where we fix \textit{any} point $y=\gamma(t)$ for a minimiser $\gamma$ defined on $[0,2t]$. By Corollary~\ref{Linfty}, this inequality yields that:
$$
|\phi(t,\gamma_1(t))-\overline{\phi}(t,\gamma_2(t))-R| \leq 2 C_3 C_4 \exp(-K_3 t),
$$
where
$$
R=\phi(t,y)-\overline{\phi}(t,y),
$$
Note that $R$ does not depend on $x$. On the other hand, using (\ref{expclose}), by Lemma~\ref{connect} we get that there is a random constant $C_7$ such that:
$$
|A(\gamma_1|_{[t,T]})-A(\gamma_2|_{[t,T]})| \leq C_7 \exp (-K_3 t).
$$
Therefore, by (\ref{first}), we get:
\begin{align} \nonumber
& |\phi(2t)-\overline{\phi}(2t)|_{L_{\infty}/\R} \leq \sup_{x \in S^1} |\phi(2t,x)-\overline{\phi}(2t,x)-R|
\\ \nonumber
& \leq (C_7+2 C_3 C_4) \exp (-K_3 t).
\end{align}
This proves the lemma's statement.\ $\qed$
\bigskip
\\ \indent
The following result follows from Lemma~\ref{expcontract} using Lemma~\ref{GN} and Lemma~\ref{W11}. Indeed, it suffices to observe that:
\begin{align*}
&|u(t)-\overline{u}(t)|_{p}=|\phi_x(t)-\overline{\phi}_x(t)|_{p}
\\ \indent
& \overset{p}{\lesssim} |\phi(t)-\overline{\phi}(t)-R|^{1/2p}_{1} |u_x(t)-\overline{u}_x(t)|^{1-1/2p}_{1}
\\ \indent
& \overset{p}{\lesssim} |\phi(t)-\overline{\phi}(t)-R|^{1/2p}_{1}.
\\ \indent
& \overset{p}{\lesssim} |\phi(t)-\overline{\phi}(t)-R|^{1/2p}_{\infty}.
\end{align*}

\begin{cor} \label{expcontractcor}
Consider two solutions $u$ and $\overline{u}$ of (\ref{Bur}) defined on the time interval $[0,+\infty)$. There exist $M>0$ such that for any $p>0$ we have:
$$
|u(t)-\overline{u}(t)|_{p} \leq C_8(p) \exp(-Mt/p),\ t \geq 0.
$$
\end{cor}
\bigskip
\indent
\textbf{Proof of Theorem~\ref{main}:} By the Fubini theorem, it suffices to prove this result in the case when the measures $\mu_1$ and $\mu_2$ are two Dirac measures concentrated on the initial conditions $u^0,\overline{u^0} \in L_p$.
\\ \indent
By a contradiction argument, it follows from Lemma~\ref{expcontract} that if we denote by $B$ the event
$$
B=\lbrace \omega \in \Omega\ |\ |u(t)-\overline{u}(t)|_{L(p)} \geq \exp(-Mt/2p) \rbrace,
$$
then  we have:
\begin{align} \nonumber
\Pe(B) & \leq \exp(-Mt/2p) \E C_8(p),\ t \geq 0.
\end{align}
Now consider a function $g$ defined on $L_p$ which satisfies $|g|_L \leq 1$. We have for $t \geq 0$:
\begin{align} \nonumber
&\E (|g(\phi(t))-g(\overline{\phi}(t))|_{p} )
\\ \nonumber
& \leq \Pe(B)\ \E (|g(\phi(t))-g(\overline{\phi}(t))|_{p}\ |\ B)
\\ \nonumber
&+ \Pe(\Omega-B)\ \E (|g(\phi(t))-g(\overline{\phi}(t))|_{p}\ |\ \Omega-B)
\\ \nonumber
& \leq 2 \Pe(B)+ \Pe(\Omega-B) \exp(-Mt/2p)
\\ \nonumber
& \leq (2 \E C_8+1) \exp(-Mt/2p).
\end{align}
This proves the theorem's statement with $C'=2 \E C_8+1$ and
\\
$K=M/2p$. \qed

\begin{rmq}
The estimate in Lemma~\ref{expcontract} is uniform with respect to the initial conditions: in other words, we have
$$
\E \sup_{\phi^0, \overline{\phi^0} \in C^0} |\phi(t)-\overline{\phi}(t)|_{L_{\infty}/\R} \leq C_6 \exp(-Mt),\ t \geq 0.
$$
A similar statement holds for the estimate in Corollary~\ref{expcontractcor}.
\end{rmq}

\section*{Acknowledgements}
\indent
I am very grateful to P.~Bernard, R.~Iturriaga, R.~Joly, K.~Khanin, S.~Kuksin, L.~Rifford, P.~Thieullen, N.~Vichery and K.~Zhang for helpful discussions. I would also like to thank separately A.~Davini, who spent a lot of time helping me to understand the paper \cite{ISM09} and A.~Fathi, who introduced we to the weak KAM community through the ANR project WKBHJ ANR-12-BS01-0020.
\bigskip

\bibliographystyle{plain}
\bibliography{Bibliogeneral}

\begin{center}
Alexandre Boritchev
\\
University of Lyon
\\
CNRS UMR 5208
\\
University Claude Bernard Lyon 1
\\
Institut Camille Jordan
\\
43 Blvd. du 11 novembre 1918
\\
69622 VILLEURBANNE CEDEX
\\
FRANCE
\\
E-mail: alexandre.boritchev@gmail.com
\end{center}

\end{document}